\documentclass[11pt,journal,onecolumn]{IEEEtran}

\usepackage[left=1.25in,right=1.25in,top=1in,bottom=1in]{geometry}
\IEEEoverridecommandlockouts                             
\usepackage{xspace,amssymb,epsfig,syntonly}
\usepackage{epsfig,amsmath,color}
\usepackage{xspace,syntonly,empheq}
\usepackage{url}
\usepackage{bbm}
\usepackage{mathtools}
\usepackage{cite}
\usepackage{graphicx}
\usepackage{algorithm,algorithmic}
\usepackage{verbatim}
\usepackage{amssymb}
\usepackage{amsthm}
\usepackage{circuitikz}
\usepackage{enumerate}
\usepackage{dsfont}
\usepackage{booktabs} 
\usepackage{epstopdf}

\newcommand{\utwi}[1]{\mbox{\boldmath $#1$}}

\renewcommand{\hat}{\widehat}
\renewcommand{\tilde}{\widetilde}

\newcommand{\cL}{{\cal{L}}}

\newcommand{\cT}{{\cal T}}

\newcommand{\cB}{{\cal B}}

\newcommand{\cX}{{\cal X}}
\newcommand{\cY}{{\cal Y}}
\newcommand{\cW}{{\cal W}}

\newcommand{\bc}{{\bf c}}

\newcommand{\bb}{{\bf b}}

\newcommand{\be}{{\bf e}}

\newcommand{\bs}{{\bf s}}
\newcommand{\bx}{{\bf x}}

\newcommand{\bv}{{\bf v}}
\newcommand{\bw}{{\bf w}}

\newcommand{\bz}{{\bf z}}
\newcommand{\by}{{\bf y}}
\newcommand{\bA}{{\bf A}}
\newcommand{\bB}{{\bf B}}
\newcommand{\bC}{{\bf C}}

\newcommand{\bJ}{{\bf J}}

\newcommand{\bL}{{\bf L}}

\newcommand{\bS}{{\bf S}}

\newcommand{\bI}{{\bf I}}
\newcommand{\bX}{{\bf X}}
\newcommand{\bZ}{{\bf Z}}
\newcommand{\bW}{{\bf W}}

\newcommand{\blambda}{{\utwi{\lambda}}}

\newcommand{\reals}{\mathbb{R}}

\newcommand{\sfT}{\textsf{T}}

\newcommand{\one}{\mathds{1}}
\DeclareMathOperator*{\argmin}{argmin}

\DeclarePairedDelimiterX{\norm}[1]{\lVert}{\rVert}{#1}

\theoremstyle{definition}

\theoremstyle{definition}

\title{\huge Optimization and Learning with Information Streams: Time-varying Algorithms and Applications} 

\author{Emiliano Dall'Anese, Andrea Simonetto, Stephen Becker, Liam Madden}

\begin{document}

\maketitle

\begin{abstract}
There is a growing cross-disciplinary effort in the broad domain of optimization and learning with streams of data, applied to settings where traditional batch optimization techniques cannot produce  solutions at time scales that match the inter-arrival times of the data points due to computational and/or communication bottlenecks.
Special types of online algorithms can handle this situation, and this article focuses on such time-varying optimization algorithms,
with  emphasis on Machine Leaning and Signal Processing, as well as data-driven Control. Approaches for the design of time-varying or online first-order optimization methods are discussed, with emphasis on algorithms that can handle errors in the gradient, as may arise when the gradient is estimated.
Insights on performance metrics and accompanying claims are provided, along with evidence of cases where algorithms that are provably  convergent in batch optimization may perform poorly in an online regime. The role of distributed computation is discussed. Illustrative numerical examples for a number of applications of broad interest are provided to convey key ideas. 
\end{abstract}

%%%%%%%%%%%%%%%%%%%%%%%%%%%%%%%%%%%
\section{Optimization and Learning With Information Streams}
\label{sec:optwithstreams}
%%%%%%%%%%%%%%%%%%%%%%%%%%%%%%%%%%%

Convex optimization underpins many important Statistical Learning, Signal Processing (SP), and Machine Learning (ML) applications. From the dawn of these fields, where least squares and kernel-based regression were prevalent across many research domains~\cite{friedman2001elements},
\footnote{Due to strict submission policies, we provide only a set of representative references.} 
to contemporary Big Data applications for online social
media, the Internet, and complex infrastructures, convex optimization enabled the development of core SP and ML algorithms and provided means to uncover   fundamental insights on  implementation trade-offs.
Modern Big Data problems, involving tasks as diverse as  kernel-based learning, sparse subspace clustering, support vector machines, and subspace tracking via low-rank models, have stimulated a rich set of research efforts that led to  breakthrough approaches for the development of scalable, efficient, and parallelizable data-processing and learning algorithms. Many of these algorithms come with a precise analysis of the convergence rates in various settings~\cite{cevher2014convex}. 

Despite these advances, the advent of streaming data sources in many engineering and science domains poses severe computational strains on existing algorithmic solutions. The ability to store, process, and leverage information from heterogeneous and high-dimensional data streams using solutions that  are grounded on batch optimization methods~\cite{NocedalWright,Be17} can no longer be taken for granted.  Timely application areas that strive for  real-time and distributed data processing and learning methods include networked autonomous vehicles,  power grids, communication systems  and the emerging Internet-of-Things (IoT) infrastructure, among many others.

In this article, we review the recent \emph{time-varying} optimization framework  \cite{popkov2005gradient,Zavala2010,Bernstein2019feedback,Paper3}, which poses a sequence of optimization problems, and  departs from batch optimization on a central processor in favor of computationally-light  algorithmic solutions where data points are processed on-the-fly and without central storage.
A first goal is to illustrate how the time-varying optimization formalism can naturally provide means to model SP and ML tasks with streams of data (with a number of instances provided in Table~\ref{tab:Applications}). The paper then addresses key questions pertaining to the design and analysis of first-order algorithms that can effectively solve these time-varying optimization problems; in particular, key aspects that will be highlighted include: i) challenges  in the design of online algorithms, along with concepts that relate the inter-arrival time of the data and the computational time of the algorithmic steps; ii) relevant metrics that can be utilized to analyze the performance of the algorithms; in particular, guidelines for the selection of given performance metrics (based on the mathematical structure of the time-varying problem) will be provided; and, iii) challenges related to distributed implementation of the online algorithms.

Before proceeding, we also point out that, beyond SP and ML, the time-varying optimization setting considered here  is relevant also for emerging data-driven control  (DDC) architectures, where learning applications are tightly-integrated components of closed-loop control systems. In this case, 
 learning components may provide means to evaluate or approximate on-the-fly constraints and costs~\cite{Bernstein2019feedback,koller2018learning}, or to drive the output of dynamical systems to solutions of time-varying optimization problems by learning first-order information of the cost functions~\cite{colombino2019online,poveda2017framework} (see Table~\ref{tab:Applications} for some examples of instances);
 this is an area that is rooted at the crossroads of SP, ML, Optimization and Control, with a natural cross-fertilization of tools and methods developed by different communities (and the divisions can be somewhat arbitrary). 
\begin{table}[h!]
  \begin{center}
    \caption{Example of instances of robust online algorithms ($^\star$: Illustrated in this article ).} 
    \vspace{-.2cm}
    \label{tab:Applications}
    \begin{tabular}{ll} 
    \toprule
    \textbf{Machine Learning / Signal Processing} & \textbf{Data-driven Closed-loop Control}\\ 
    \midrule
     Subspace  tracking$^\star$, robust subspace  tracking~\cite{dixit2019online}  & Measurement-based online algorithms$^\star$~\cite{Bernstein2019feedback} \\
     Subspace clustering, sparse subspace clustering$^\star$~\cite{friedman2001elements,elhamifar2013sparse} & Extremum seeking~\cite{poveda2017framework} \\
    Sparse~\cite{cevher2014convex}, kernel-based~\cite{friedman2001elements}, robust, linear regression~\cite{friedman2001elements} & Online optimization as feedback controller~\cite{colombino2019online} \\
    Zeroth-order methods~\cite{Hajinezhad19},  bandit methods~\cite{chen2018bandit,besbes2015non}  &  Predictive control with Gaussian Processes~\cite{koller2018learning}   \\
    Support vector machines~\cite{friedman2001elements} & \\
    Learning problems over networks~\cite{onlineSaddle} & \\
    \bottomrule 
    \end{tabular}
    \vspace{-.4cm}
  \end{center}
\end{table}

\textbf{Batch/convergence vs online/tracking.} A central  question related to ML and DDC applications pertains to the  ability of existing iterative optimization algorithms --- especially first-order methods~\cite{cevher2014convex,Be17} --- to  handle data streams effectively. 
Suppose that data points arrive sequentially at intervals of $\delta > 0$ seconds; discretize the temporal axis as $t \in \cT := \{ k \delta: \delta > 0, k \in \mathbb{N} \}$, where $\delta$ can be selected as the inter-arrival time of data, and suppose  that a given ML task involves the solution of the following  problem at time $t \in \cT$, based on data $\bZ_{t} = \{\bz_\tau, \tau \in \cW_t\}$ gathered over a (possibly sliding) window $\cW_t \subset \cT$:
\vspace{-.2cm}
\begin{align}
\label{eq:time_varying_problem}
\text{Find}\; f_t^* := \min_{\bx_t \in \cX_t} f_t(\bx_t; \bZ_{t})  
\quad\text{or}\quad
\bx_t^* \in \argmin_{\bx_t \in \cX_t} f_t(\bx_t; \bZ_{t})  
\end{align}
where $\bx_t^* \in \mathbb{R}^n$ is the parameter of interest,
$\cX_t \subseteq \mathbb{R}^n$ is a
non-empty closed convex set, and the time-varying function may take the form $f_t(\bx_t; \bZ_{t})  := h_t(\bx_t; \bZ_{t}) + g_t(\bx_t)$, with $h_t$ convex and $L_t$-smooth (i.e, $\nabla h_t$ is $L_t$-Lipschitz continuous) and $g_t$ convex but not necessarily smooth. For example, to illustrate the temporal variability of the function as well as the explicit dependency of the cost on  the data stream,  for an $\ell_1$-regularized least-squared problem one may have $f_t(\bx_t; \bZ_{t}) = \sum_{\tau \in \cW_t} \|\bA_\tau \bx_t - \bb_\tau\|_2^2 + \lambda_t \|\bx_t\|_1$, with $\lambda_t > 0$ a time-varying sparsity-promoting parameter, $\bZ_\tau = \{\bA_\tau, \bb_\tau, \tau \in \cW_t\}$, and $\cW_t = \{t - W, \ldots, t\}$. As another example, problem~\eqref{eq:time_varying_problem} could be used for subspace tracking based on a sliding window of (vectorized) video images, by setting $h_t$ to be a least-squares term, and $g_t$ a nuclear norm regularization~\cite{dixit2019online} (these examples will be illustrated shortly). Hereafter, to simplify the notation, we drop $\bZ_t$ from the cost function, letting the dependency of $f_t$ on the data be implicit.

Supposing that~\eqref{eq:time_varying_problem} is solved using a proximal gradient method or an accelerated proximal gradient method (proximal gradient methods are generalizations of projected gradient methods, see \cite{cevher2014convex,Be17}), it is known that when 
$f_t$ is convex and $h_t$ is $L$-smooth,
the number of iterations required to obtain an objective function within an error $\epsilon$ is $\mathcal{O}(L \|\bx_{t,0} - \bx_t^*\|^2 / \epsilon )$ and $\mathcal{O}(\sqrt{L \|\bx_{t,0} - \bx_t^*\|^2 / \epsilon} )$,
respectively, with $\bx_{t,0}$ the starting point for the algorithm, and $\bx_t^*$ the optimal solution at time $t$~\cite{cevher2014convex,Be17} ($\mathcal{O}$ refers to the big O notation).
If one can computationally afford a number of proximal gradient steps in the above order within an interval $\delta$, then it is clear that $\bx_t^*$ can be identified, within an acceptable error, at each step $t$. However, what if data points arrive at a rate such that only a few steps (or even just one step)  can be performed before a new datum arrives?   

Continuing to use the proximal gradient method as a running example, and taking the extreme (yet realistic) case where only one step can be performed within an interval $\delta$, conventional wisdom would suggest to utilize the following \emph{online} implementation: 
\begin{align}
\bx_{t} = \textrm{prox}_{\alpha_t g_t, \cX_t}\{\bx_{t-1} - \alpha_t \nabla_{\bx} h_t(\bx_{t-1}) \}, \hspace{.5cm} t = 1, 2,  \ldots 
\label{eq:forward-backward}
\end{align}
where $\textrm{prox}_{\alpha g, \cX}\{\by\} := \arg \min_{\bx \in \cX} \left\{ g(\bx) + \frac{1}{2 \alpha} \|\bx - \by\|_2^2 \right\}$ is the proximal operator~\cite{Be17},  $\{\alpha_t\}$ is a given step-size sequence, and $f_t$ is built based on new data up to time $t$ (and possibly discarding older data, as in a sliding window mode)\footnote{We stress that it may be possible to perform multiple proximal-gradient steps within an interval $\delta$, but we consider the case of one step to simplify the notation.}.  Relevant questions in this setting revolve around the  \emph{definition of suitable metrics} to  analytically characterize the performance of the algorithm~\eqref{eq:forward-backward},
since the classical notions of convergence and $\epsilon$-accurate solutions for batch optimization are no longer suitable.
Viewing~\eqref{eq:time_varying_problem} under the lens of a time-varying optimization formalism~\cite{Zavala2010,Paper3,Bernstein2019feedback,onlineSaddle,dixit2019online}, metrics that are related to \emph{tracking} of the sequence  $\{f_t^*, t \in \mathbb{N}\}$ or 
sequences of optimal solutions $\{\bx_t^*, t \in \mathbb{N}\}$ will be discussed in Section~\ref{sec:performance}. 

\begin{figure}[t!]
    \centering
    \includegraphics[height=4.8cm]{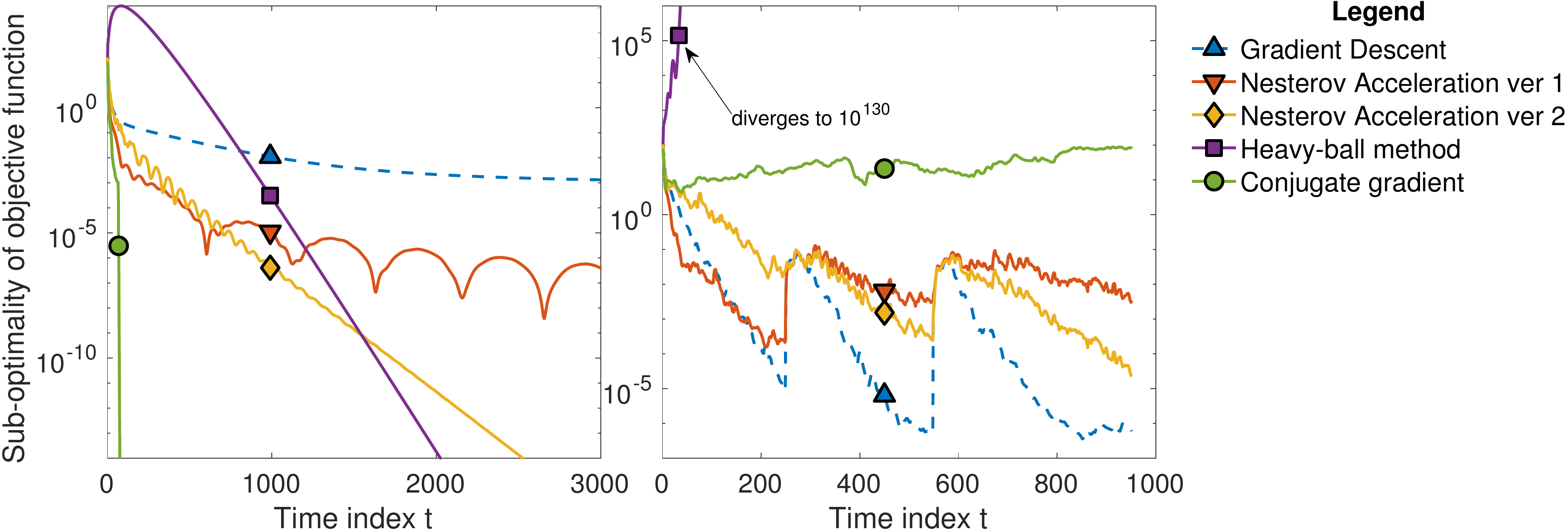}
    \vspace{-.3cm}
        \caption{Example of a 50 dimensional time-varying least-squares problem, defined using a sliding window of 50 data points, for 950 time points; two big jumps in the solution near time indices 250 and 550 (by design).
    Left: Convergence in the \emph{static} case; Right: plot shows the performance of various algorithms on \emph{tracking} the optimal objective value.
    Nesterov ver.~1 does not use knowledge of strong convexity, while ver.~2 does. The non-linear conjugate gradient exploits the quadratic objective to have an exact line-search (usually impractical), and is the variant from 
    \cite[Eq.~(5.49)]{NocedalWright}. }
    \label{fig:dynamic}
\end{figure}

Before proceeding, 
we bring
up a point that highlights another challenge in designing and analyzing online algorithms. One may surmise that a na\"{i}ve online implementation of algorithms conceived for batch computation may just work well,
with algorithms that are faster for batch optimization still being faster in time-varying optimization.
However, this is not always the case. Surprisingly, the best algorithms in the static case may be the worst algorithms in the dynamic case, as shown in our illustrative numerical results in Fig.~\ref{fig:dynamic}. The 
heavy ball method 
can even diverge for a simple time-varying least-squares problem. 

%%%%%%%%
%%%%%%%%

\textbf{Time-varying optimization vis-\`{a}-vis online learning.} The time-varying optimization formalism 
is closely aligned with existing works on  ``online learning in dynamic environments''   ~\cite{dixit2019online,Jadbabaie2015,hall2015online,shalev2012online,hazan2016introduction}  from a basic mathematical standpoint. 
However, a key conceptual difference is that the online  algorithms for  time-varying optimization described here are ``\emph{computation limited},'' whereas online learning  is ``\emph{data limited}'' or ``\emph{information limited}'' (but not necessarily  computation limited).  To understand the term ``\emph{data limited},'' take as an example the online learning framework in a learner-environment setting outlined in~\cite{Jadbabaie2015,hall2015online}. In this case, the online learning counterpart of~\eqref{eq:forward-backward} could be written as $\bx_{t} = \textrm{prox}_{\alpha_t g_{t-1}, \cX} \{\bx_{t-1} - \alpha_t \nabla_{\bx} h_{t-1}(\bx_{t-1}) \}$; that is, a typical online learning algorithm produces a ``prediction'' $\bx_{t}$ based on information of the cost function (in the form of functions or functional evaluations) and data  that are available up to the time $t-1$ (although there is no consensus on the terminology ``learning'' vs ``optimization'' in the literature for this learner-environment  setting, we will hereafter use the term ``learning'' for algorithms that predict $\bx_{t}$ based on $f_{t-1}$; we will use the term ``optimization'' for the case where $f_{t}$ is known at time $t$). Once the prediction $\bx_{t}$ is computed, partial or full feedback regarding the function $f_{t}$ is revealed. The performance of online learning is therefore evaluated relative to the best action in hindsight; that is, relative to the case where $f_{t}$ is available to the learner. Furthermore, the majority of the online learning frameworks assume that the set $\cX_t$ is static; this is to avoid  infeasibility of $\bx_{t}$ that one may have if $\cX_{t}$ is unknown to the learner at time $t$.   

On the other hand, the time-varying optimization setting outlined here is mainly driven by computational bottlenecks. At time $t$, information regarding the function $f_t$ to be minimized is available in terms of either its functional form or (an estimate of) its gradients. The algorithm then seeks to obtain an optimal solution $\{\bx_t^*\}$; however, because of complexity and data rate considerations,  only one or a few algorithmic steps can be performed within $\delta$ seconds, before a new datum $\bz_{t+1}$ arrives  (and the process then restarts in order to seek a new optimizer $\{\bx_{t+1}^*\}$). While we allow for multiple algorithmic steps within $\delta$ seconds, online learning methods such as the online mirror descent~\cite{hall2015online,Jadbabaie2015} typically consider only one step. As explained shortly in Section~\ref{sec:performance}, the performance of a time-varying algorithm  is measured against the  solution that would have been obtained had we had the time to run an algorithm to convergence at each interval $\delta$. In this time-varying setting, concepts that relate the inter-arrival time   of the data, the sampling time $\delta$, and the computational time of the algorithms play a key role, and clearly assumptions must be made about how $f_t$ varies over time. For example, there might be enough time to compute two gradients of $f_t$, or one gradient and two function evaluations (for a line search), so an algorithm can choose how to spend  computational resources.

\begin{figure}[t]
    \centering
    \includegraphics[width=3.8in]{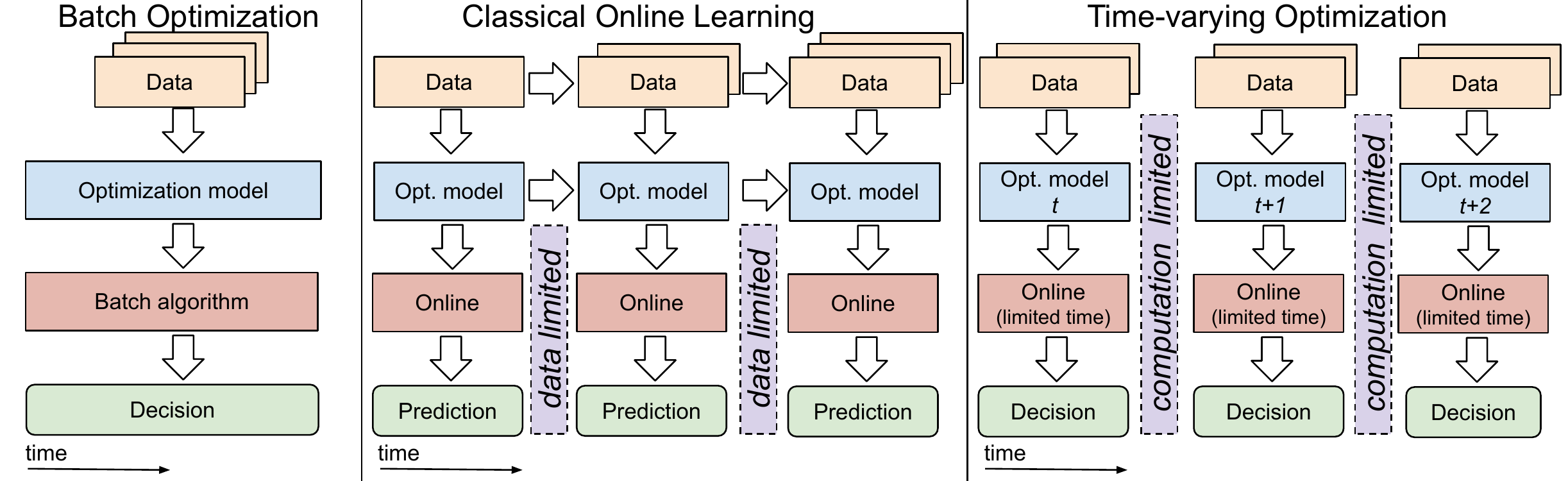}
    \includegraphics[width=2.0in]{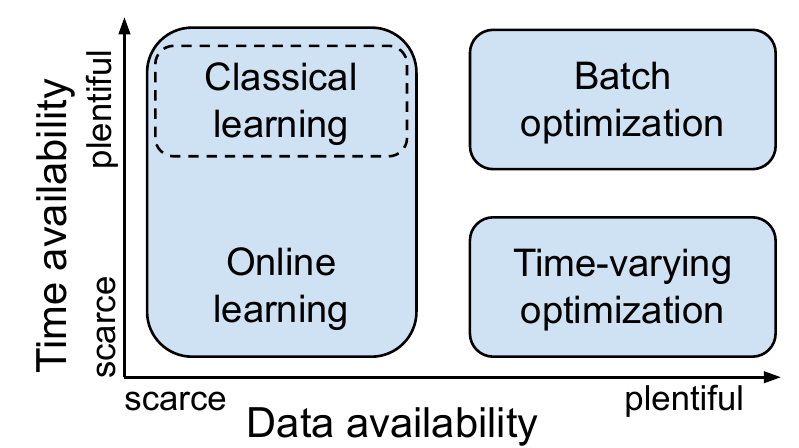} 
    \vspace{-.2cm} 
    \caption{Online learning and time-varying optimization are both sequential, but in online learning data is limited in the sense that $f_t$ is not available. In time-varying optimization, there is no restriction of information, but the full problem cannot be solved within a single time step due to computational cost. The right-most plot provides an illustrative distinction between: (i) time-varying algorithms, that have knowledge of $f_t$, but are computationally-limited; (ii) classical learning methods (such as the follow-the-leader method) that do not necessarily have computational limitations~\cite{shalev2012online,hazan2016introduction}; and, (iii) online learning, where only one algorithmic step is performed per time interval (such as the online mirror descent method~\cite{hall2015online}). In both (ii) and (iii), no knowledge of $f_{t}$ is available.   
    }
    \label{fig:my_label}
\end{figure}

It is also worth mentioning that, while the considered online learning example involved a proximal-gradient step (to highlight the subtleties relative to~\eqref{eq:forward-backward}),  computationally-heavier learning algorithms can be utilized to produce the prediction $\bx_t$, with  the computational effort between time-slots not necessarily being a  concern~\cite{shalev2012online,hazan2016introduction}. For example, the popular follow-the-leader method requires a batch solution of an optimization problem at each time step if the cost functions are not linear~\cite{shalev2012online}. In Fig.~\ref{fig:my_label}, we refer to this setting as ``classical learning''.

We focus here on algorithms implemented with a  constant step-size; this is a natural choice for cases where the optimal solution $\{\bx_t^*\}$ remains transient and the  algorithm  runs indefinitely. This is another  distinction relative to online learning with a time-invariant optimizer, where the  step-size may depend on the time-horizon or a ``doubling trick''~\cite[Sec. 2.3.1]{shalev2012online} is utilized (with the latter still involving changes in the algorithm based  on how many iterations have been taken). 
An example to distinguish standard online learning from time-varying optimization is spam filtering.  At time $t$, a user receives an email message with features $\bz_t$, and their email software must decide whether to label the email as a spam message or a legitimate message.  An \emph{online learning} problem in this scenario is to make the best use of prior emails (and their correct labels) to make a prediction for the new email, after which the user will supply the correct answer, and the software will take this into account at time $t+1$. A \emph{time-varying} problem in this scenario is the case where we assume that all users receive the same type of spam email and thus do not need an individually trained classifier, but that the nature of email spam evolves over time, and hence the email  provider must update the spam classifier every day. The email provider has access to all the emails of their users, and thus plenty of data, but might use a complicated classifier that cannot be completely trained in one day.

%%%%%%%%%%%%%%%%%%%%%%%%%%%%%%%%%%%
\section{Time-Varying Algorithms}
\label{sec:algorithms}
%%%%%%%%%%%%%%%%%%%%%%%%%%%%%%%%%%%

We will overview online algorithms to track solutions of  time-varying problems of the form~\eqref{eq:time_varying_problem} that are designed based on three key principles:  

\noindent \emph{[P1] First-order methods}. First-order methods are computationally-light, they facilitate the derivation of parallel and distributed architectures, and they can handle non-smooth objectives by leveraging the proximal mapping~\cite{cevher2014convex,Be17}. In the context of this article, they exhibit  robustness to inaccuracies in the gradient information --- an important feature further explained next\footnote{In some scenarios, 
second-order information is computationally cheap to obtain and then second-order methods are competitive (cf.~\cite{Paper3} in the context of prediction-correction methods) but we do not pursue these methods here.}.

\noindent \emph{[P2] Approximate first-order information}. We consider first-order algorithms 
that are robust to inaccurate gradient estimates;
more precisely, the online algorithm is executed using a sequence  $\{\bv_t, t \in \mathbb{N}\}$, with well-posed assumptions on the sequence of differences  $\be_t :=  \nabla_{\bx} h_t(\bx_{t-1}) - \bv_t$. For example,  assumptions relative to $\be_t$ in existing literature involve boundedness of $\|\be_t\|$ (in a given norm)~\cite{Dallanese2018feedback,schmidt2011convergence}, as explained in Section~\ref{sec:performance}. In a stochastic case, boundedness of $\mathbb{E}[\|\be_t\|]$ (where $\mathbb{E}[\cdot]$ denotes expectation) is presumed~\cite{dixit2019online}. This setting finds important applications in ML and DDC with data streaming, with prime examples outlined shortly. 

\noindent \emph{[P3] Distributed computation}. We cover problems of the form~\eqref{eq:time_varying_problem} or suitable reformulations that are to be solved collaboratively by a network. 

We start by revisiting the proximal gradient method~\eqref{eq:forward-backward} under the lens of \emph{[P1]--[P2]}. This algorithm is relevant for a number of  instances listed in Table~\ref{tab:Applications}, in particular when $\cX_t = \mathbb{R}^n$ or projection onto the sets $\{\cX_t, t \in \cT\}$ is computationally cheap. We then turn the attention to primal-dual methods and variants~\cite{Koshal11}; these methods  are naturally applicable to the case where 
the set $\cX_t$  in~\eqref{eq:time_varying_problem} is expressed as $\cX_t = \cY_t \cap \{\bx: \bc_t(\bx) \leq \mathbf{0}\}$, with $\cY_t$ convex (and involving a computationally-cheap projection), and $\bc_t$ a vector-valued convex function; here, the constraint $\bc_t(\bx) \leq \mathbf{0}$ is  dualized to construct the Lagrangian function. This setting is relevant, for example, in network optimization problems with data streams~\cite{onlineSaddle,Koshal11,Bernstein2019feedback}. On the other hand, a similar structure emerge in consensus-based reformulations of~\eqref{eq:time_varying_problem}, where $\cX_t = \cY_t \cap \{\bx: \bC \bx = \mathbf{0} \}$, with $\bC$ a given consensus matrix constructed based on the  communication graph~\cite{Ling14admm,dimakis2010gossip,boyd2011distributed}.

\smallskip

\subsection{Proximal gradient method}
\label{sec:proximal}

The online algorithm~\eqref{eq:forward-backward} with approximate first-order information amounts to the sequential execution of the following step:
\begin{align}
\bx_{t} = \textrm{prox}_{\alpha_t g_t, \cX_t}\{\bx_{t-1} - \alpha_t \bv_t \}, \hspace{.5cm} t = 1, 2,  \ldots 
\label{eq:forward-backward-v}
\end{align}
\noindent where we 
recall that 
$\bv_t \approx \nabla h_t(\bx_t)$.
If $g_t=0$, \eqref{eq:forward-backward-v} reduces to the online projected gradient method with approximate gradient information.
We focus our attention to algorithms with a constant step-size; that is, $\alpha_t = \alpha > 0$ for all $t \in \mathbb{N}$; this is reasonable when no prior on the evolution of the problem is available, and the algorithm is executed indefinitely (as opposed to a given finite interval).
The algorithm \eqref{eq:forward-backward-v} is the starting point for all variants we consider, and handles the key issues of the time-varying setting, namely 
functions and constraints (and hence the optimizers) possibly changing at each step $t$, and inexact gradients.
The algorithm \eqref{eq:forward-backward-v} is utilized in the examples of applications presented  in Section~\ref{sec:Examples_applications}. 
Section~\ref{sec:Handling_Lagrangian_functions} will discuss a more general algorithm based on Lagrangian functions, and Section~\ref{sec:performance} will elaborate on performance metrics.

\subsection{Examples of applications of the proximal gradient method}
\label{sec:Examples_applications}

\noindent $\bullet$ {\bf \emph{ML example \#1: Subspace tracking for video streaming}}. 
Robust principal component analysis (PCA) can be used to separate foreground from background in video, among many other applications. The model is that a matrix $\bZ\in\reals^{p^2 \times T}$, which encodes a video as $p^2$ pixels by $T$ video-frames, can be decomposed as $\bZ \approx \bS + \bL$ where $\bS$ is sparse (foreground) and $\bL$ is low-rank (background). Now suppose $\bZ_t$ is a video clip, and $\bZ_{t+1}$ is the subsequent video clip, and the objective is to decompose \emph{all} video clips into foreground and background in a streaming and \emph{real-time} fashion. This form of \emph{robust subspace tracking} was considered by \cite{dixit2019online}, and is modeled by solving the following problem (for parameters $\lambda,\rho > 0$):
\begin{equation} \label{eq:RPCA}
    \min_{\bL_t,\bS_t}\; \lambda\|\bS_t\|_1 + \|\bL_t\|_* + \frac{\rho}{2}\|\bL_t + \bS_t - \bZ_t\|_F^2
    =
    \min_{\bS_t}\; \underbrace{\lambda\|\bS_t\|_1}_{g_t(\bS_t)}
    + \underbrace{\min_{\bL_t}\; \|\bL_t\|_* +
    \frac{\rho}{2}\|\bL_t + \bS_t - \bZ_t \|_F^2}_{h_t(\bS_t)}
\end{equation}

\begin{figure}[t]
    \centering
    \includegraphics[width=2.7in]{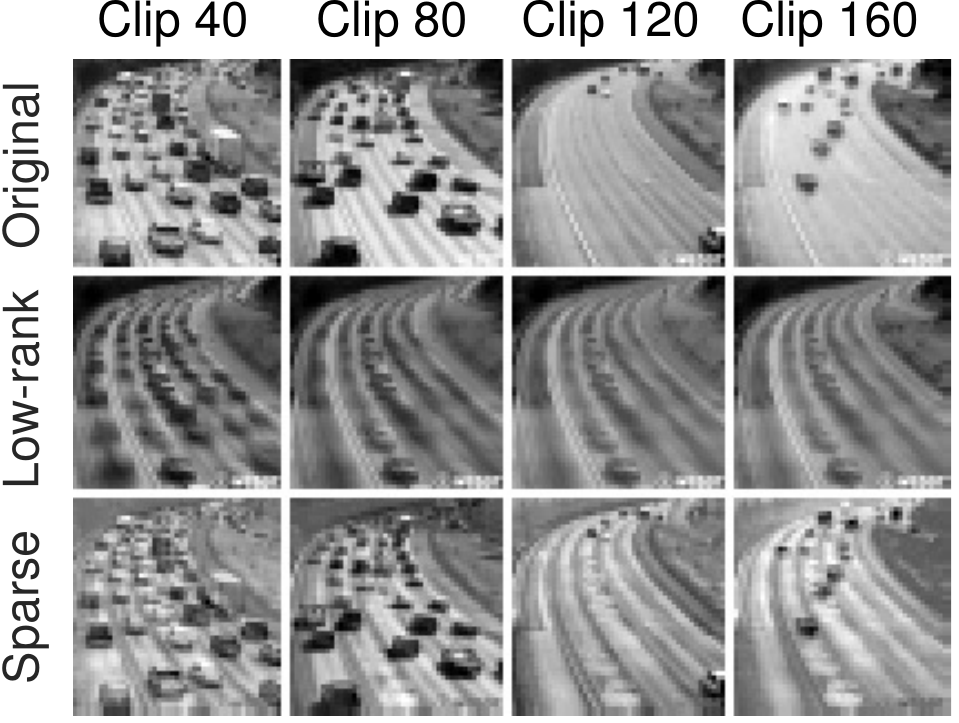}
    \includegraphics[width=2.7in]{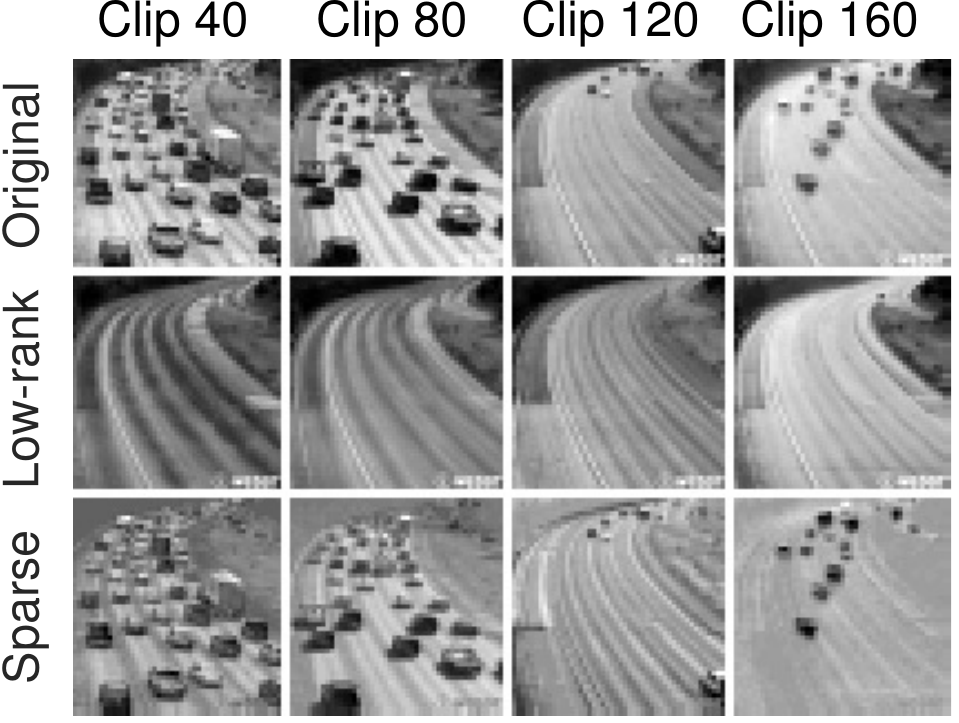}
    \caption{Results of robust PCA on traffic footage. Left: 2 iterations of proximal gradient per time video clip. Right: 10 iterations of proximal gradient per video clip. \vspace{-.4cm}}
    \label{fig:traffic}
\end{figure}

\noindent where the term $\|\bL_t\|_*$ is the nuclear norm. 

Identifying $g_t$ and $h_t$ as in Eq.~\eqref{eq:RPCA}, 
then $h_t$ is  
differentiable and $\nabla h_t$ is $\rho$-Lipschitz continuous~\cite[Prop.\ 12.30]{CombettesBook2}, and 
$\nabla h_t(\bS_t) = \rho( \bL + \bS_t - \bZ_t)$
where $\bL$ solves the inner minimization problem and can be computed using the singular value decomposition (SVD) of $\tilde{\bZ} := \bZ_t-\bS_t$~\cite[Ex.~24.69]{CombettesBook2}. To speed up computation, the SVD algorithm may be allowed to produce small errors, such as in randomized SVD methods~\cite{halko2011finding}, iterative methods like Lanczos, or in methods with large roundoff error. In particular, if $T\ll p^2$, an efficient SVD algorithm is to compute the eigenvalue decomposition (EVD) of $\tilde{\bZ}^T\tilde{\bZ}$; this multiplication and EVD operation have the same asymptotic flop count as the usual SVD, but is faster since it has smaller constants and can exploit well-tuned matrix multiply algorithms (especially true on the GPU). It has higher numerical error due to squaring the condition number of $\tilde{\bZ}$. There are other choices for defining $g_t$ and $h_t$, such as those in \cite{dixit2019online}, but our choices fit into the approximate gradient framework which comes with guarantees. Eq.~\eqref{eq:RPCA} fits into Eq.~\eqref{eq:time_varying_problem} by letting $\bx_t$ be a vectorized version of $\bS_t$ with $\cX=\reals^n$ and $n=p^2T$, and~$f_t = g_t + h_t$.

To illustrate the example, we took a dataset of 254 traffic video clips from
\cite{trafficDBpaper2005}, of resolution $48 \times 48$ so $p^2=2304$, and $T$ in the range of 48 to 52, and chose $\lambda=1/p$ and $\rho=10$. Most video clips are from subsequent times, but there is a jump between the last clip in the evening and the first clip in the morning. 
On a 2-core laptop, a single proximal gradient step takes about $0.0022$ seconds. Thus for real-time processing, assuming a frame rate of 25 frames/second, one could take just under 20 iterations per video clip. Fig.~\ref{fig:traffic} illustrates the proximal gradient algorithm taking just 2 iterations per clip (left) and 10 iterations per clip (right). Obviously taking more iterations per video clip leads to better results. The quality of the background component is much better for clip 160 than it is for clip 40, also as expected.

\noindent $\bullet$ {\bf \emph{ML example \#2: Online sparse subspace clustering.} } Subspace tracking  identifies a shared low-dimensional subspace that explains most of the data~\cite{friedman2001elements}; on the other hand,  subspace \emph{clustering} is a key ML application utilized to group data points that lie in low-dimensional affine spaces~\cite{friedman2001elements}. Subspace clustering is useful when data points lie near low-dimensional affine spaces and represent qualitatively different real-world objects based on their respective affine spaces. 

\begin{figure}[H]
    \centering
    \includegraphics[width=.90\linewidth]{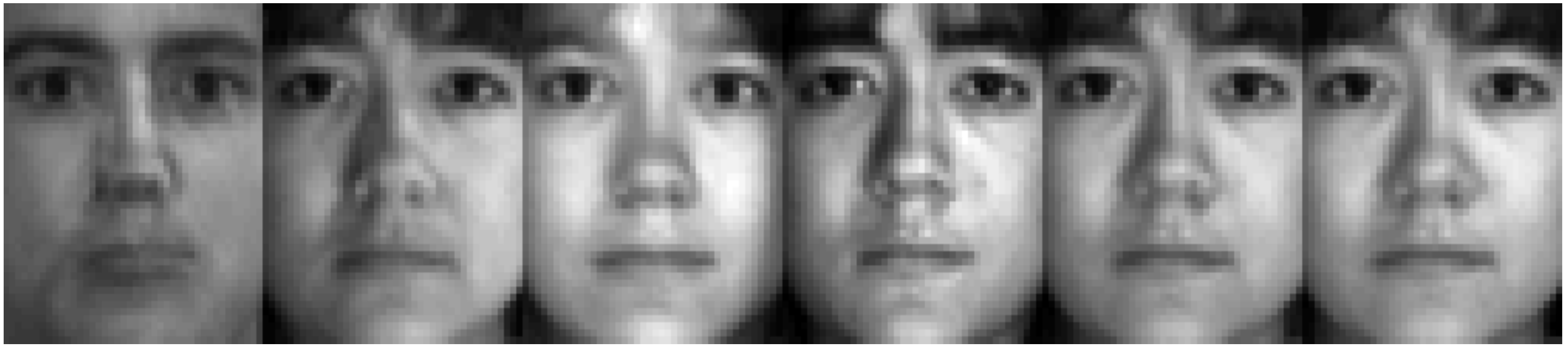}
    \vspace{-.4cm}
    \caption{
    The average image of data points labeled as representing a certain person as ``time'' increases and convergence to tracking error is reached, for times $\{1,40, 80, 120, 160, 200\}$.  ``Time'' refers to which sliding window of the data set is being used.
    }
    \label{fig:faces}
\end{figure}
Here we illustrate an approach referred to as sparse subspace clustering~\cite{elhamifar2013sparse}; it involves the solution of a sparse representation problem, followed by spectral clustering applied to the  graph corresponding to the similarity matrix formed using the minimizer of the sparse representation problem. In particular, the sparse representation problem~\cite{elhamifar2013sparse} is
\begin{align}
\label{eq:sparse-rep}
    \min_{\bX_t \in \reals^{N\times N}}\; \underbrace{\lambda\|\bX_t\|_1}_{g_t(\bX_t)}+\underbrace{\frac{1}{2}\|\bZ_t\bX_t-\bZ_t\|_F^2}_{h_t(\bX_t)}\text{ s.t.} \underbrace{\text{diag}(\bX_t)={\bf 0}\text{ and }\bX_t^T\one=\one}_{\cX_t}
\end{align}
where $\|\cdot\|_1$ is the vector $\ell_1$-norm, diag$(\bX_t)$ is the vector consisting of the diagonal elements of $\bX_t$, and ${\bf 0}$ and $\one$ are vectors of all zeros and ones respectively. 
Letting $\bx_t$ be a vectorized version of $\bX_t$ with $n=N^2$, this fits in the framework of~\eqref{eq:time_varying_problem}.
$\bZ_t$ is a sliding window so that the clustering problem does not grow in the number of data points that need to be labeled, in order to avoid creating a growing computational demand.
Fig.~\ref{fig:faces} visually represents the evolution of the center of one subspace, by averaging the data points of one cluster and showing how it changes over time, starting with a mixture of faces and then  converging on the identity of a single person.

%%%%%%%%
\begin{figure}[t]
    \centering
    \includegraphics[height=5.3cm]{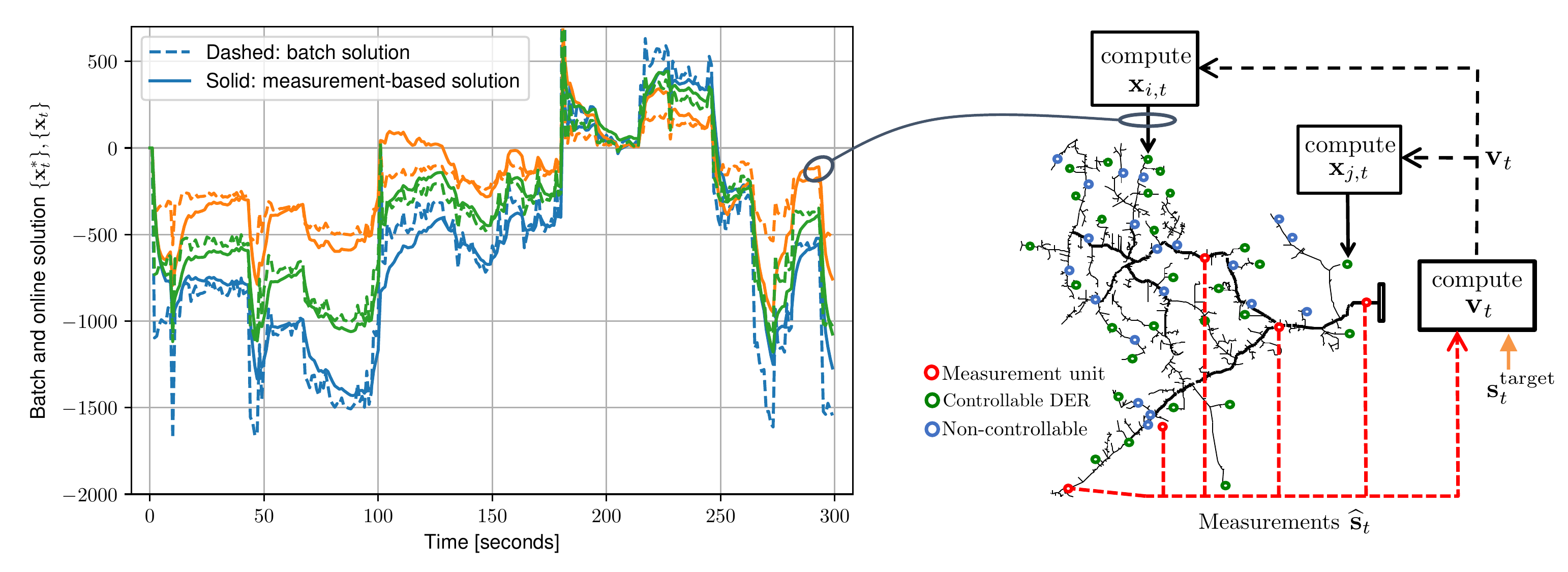}
    \vspace{-.5cm}
        \caption{Example of measurement-based online algorithm, with application to power systems (adapted from~\cite{Dallanese2018feedback}). The numerical results correspond to the case where $\bx_t = [\bx_{1,t}^\sfT, \ldots, \bx_{N,t}^\sfT]^\sfT$ are power commands for $N$ distributed energy resources (DERs), $\bw_t$ collects the powers consumed by non-controllable loads, and $s_t$ is in this case a scalar representing the setpoint for the aggregate power of both DERs and loads; i.e., $s_t = \mathbf{1}^\sfT \bx_t + \mathbf{1}^\sfT \bw_t$. Left: example of tracking performance for the real power of three representative DERs (in W). Right:  with $\bv_t$ computed by a central network operator, the update~\eqref{eq:forward-backward-v} decouples across nodes when $g_t(\bx_t) = \sum_{i = 1}^N g_{i,t}(\bx_{i,t})$ (with $N$ the number of DERs, and $\cX_t = \cX_{1,t} \times \cX_{2,t} \times \ldots \times \cX_{N,t}$) .  } 
    \label{fig:F_tracking_power_systems}
\end{figure}
%%%%%%%%

\noindent $\bullet$ {\bf \emph{DDC example \#1: Measurement-based online network optimization}}. Consider a physical network (e.g., a power system, a transportation network, or a communication network) described by a input-output map $\bs_t = \bA \bx_t + \bB \bw_t$, with $\bx_t$ a vector of control inputs, $\bw_t \in \mathbb{R}^{w}$ (possibly unknown) exogenous inputs, and  $\bA \in \mathbb{R}^{s \times n}$, $\bB \in \mathbb{R}^{s \times w}$ known network matrices; the vector $\bs_t$ collects network outputs. For example, in power systems, $\bw_t$ collects the power consumed by uncontrollable loads throughout the (possibly very large) network, $\bx_t$ collects controllable power injections, $\bs_t$ collects power flows, and the network map is based on a linearized power flow model~\cite{Dallanese2018feedback}. As an illustrative example, suppose that the function $h_t$ in~\eqref{eq:time_varying_problem} is $h_t(\bx) = \frac{1}{2}\|\bA \bx_t + \bB \bw_t - \bs_t^{\textrm{target}}\|_2^2$ in an effort to drive the network output towards a time-varying reference point $\bs_t^{\textrm{target}}$. The gradient of $h_t$ in this case reads $\nabla_{\bx} h_t(\bx_t) = \bA^\sfT (\bA \bx_t + \bB \bw_t - \bs_t^{\textrm{target}})$. Evaluating the gradient requires one to measure or estimate the vector $\bw_t$ at each step of the algorithm, and this task can be problematic (if not even impossible) in many real-world applications. If, on the other hand, sensors are deployed to measure the network output $\bs_t$, then a proxy of the gradient 
can be constructed as $\bv_t = \bA^\sfT (\hat{\bs}_t - \bs_t^{\textrm{target}})$, with $\hat{\bs}_t$ a measurement of $\bs_t$~\cite{Bernstein2019feedback}. Because of the inherent measurement errors, but also because of a possibly inaccurate knowledge of the model matrix $\bA$,
$\bv_t$ is a noisy version of $\nabla_{\bx} h_t(\bx_t)$. An example of measurement-based architectural framework is illustrated in Figure~\ref{fig:F_tracking_power_systems}, for an application to power systems (adapted from~\cite{Dallanese2018feedback}).

\subsection{Handling Lagrangian functions}
\label{sec:Handling_Lagrangian_functions}
Consider now the case where $\cX_t = \cY_t \cap \{\bx: \bc_t(\bx) \leq \mathbf{0}\}$; we recall that his setting is relevant, for example, in network optimization problems with data streams~\cite{Koshal11,Bernstein2019feedback} or in SP and ML applications where the projection onto $\cX_t$ can be computationally intensive. Focusing first on the case where $f_t(\bx_t)$ is strongly convex, we consider the design of first-order algorithms based on the \emph{time-varying} saddle-point problem~\cite{Koshal11}:  
\begin{align}
\min_{\bx_t \in \cY_t} \max_{\blambda_t \in \mathbb{R}_+^m} \cL_t(\bx_t,\blambda_t) := f_t(\bx_t) + \blambda_t^\sfT \bc_t(\bx_t) - \frac{r}{2} \|\blambda_t\|_2^2
\label{eq:reg-lagrangian}
\end{align}
where $\cL_t(\bx_t,\blambda_t)$ is a regularized Lagrangian function, $r \geq 0$ is a regularization parameter, and $\blambda_{t} \in \mathbb{R}_+^m$ is the vector of  multipliers associated with the constraint $\bc_t(\bx) \leq \mathbf{0}$. Accordingly, based on the principles \emph{[P1]--[P2]}, an approximate online primal-dual method is of the form~\cite{Bernstein2019feedback}: 
\begin{align}
\label{eq:primal-dual}
\bx_{t}  = \textrm{prox}_{\alpha g_t, \cY_t}\{\bx_{t-1} - \alpha \bv^{(L)}_t \}, \quad
\blambda_t  = \textrm{proj}_{\mathbb{R}^m_+} \left\{ (1-r) \blambda_{t-1}  + \alpha \bv^{(c)}_t  \right\}
\end{align}

\vspace{-.1cm}

\noindent where $\bv^{(c)}_t$ is a proxy for the gradient of $\cL_t$ with respect to $\blambda_t$, and $\bv^{(L)}_t$ is an estimate of $\nabla_\bx h_t(\bx) + \bJ_{\bc_t}^{\sfT}(\bx)\blambda_{t-1}$, with $\bJ_{\bc_t}$ denoting the Jacobian of $\bc_t$ (i.e., an estimate of the gradient of the smooth part of the regularized Lagrangian). It is important to notice that if $r = 0$, then~\eqref{eq:reg-lagrangian} reverts to the standard Lagrangian function; then,~\eqref{eq:primal-dual} can be utilized to track optimal primal-dual \emph{trajectories} of~\eqref{eq:time_varying_problem} based on metrics grounded on the notion of regret, but there are no linear convergence results due to the lack of strong convexity of the dual problem~\cite{Koshal11,Du2019}.
When $r > 0$, then $\cL_t(\bx_t,\blambda_t)$ becomes ${r}$-strongly concave in $\blambda_t$ and linear convergence results are available at the cost of tracking an approximate solution. This will be discussed in Section~\ref{sec:performance}.

%%%%%%%%%%%%%%%%%%%%%%%%%%%%%%%%%%%
\section{Performance Analysis: Which Metrics?}
\label{sec:performance}
%%%%%%%%%%%%%%%%%%%%%%%%%%%%%%%%%%%

Given the temporal variability of the underlying optimization problems, the so-called \emph{dynamic regret}~\cite{besbes2015non,hall2015online,onlineSaddle,Jadbabaie2015,dixit2019online,li2018using}, and (a slightly modified notion of) $Q$-\emph{linear convergence}~\cite{Paper3,Bernstein2019feedback,dixit2019online} are metrics that can be used to characterize the performance of time-varying algorithms.  This section discusses these metrics and relevant bounding techniques. Further, this section provides guidelines on the adoption of given performance metrics, based on the mathematical structure of the time-varying problem.

It is first necessary to make assumptions on
(i) a ``measure'' of the temporal variability of~\eqref{eq:time_varying_problem}, and (ii) the ``correctness'' of the first-order information $\{\bv_t\}$. 
Measures of the latter include \cite{schmidt2011convergence}:
\begin{align}
\label{eq:timevarying_sequence1}
 E_t := \sum_{\tau = 1}^{t} e_\tau \quad \text{and/or}\quad 
 E_t^\prime := \sum_{\tau = 1}^{t} e_\tau^2, \quad \text{where}\quad 
 e_t := \|\nabla_{\bx} h_t(\bx_{t-1}) - \bv_t\|,
\end{align}
with  $E_t$ and $E_t^\prime$ representing the error accumulated up to time $t$. Stochastic counterparts of~\eqref{eq:timevarying_sequence1} of the form $\mathbb{E}[e_t]$ are also considered as discussed in, e.g.,~\cite{dixit2019online}.

Temporal variability of the problem~\eqref{eq:time_varying_problem} could be measured based on how fast its optimal solutions 
evolve; more precisely, assuming first that the cost in~\eqref{eq:time_varying_problem} is strongly convex at all times
(and, therefore, the trajectory $\{\bx_t^*, t\in \cT\}$ is unique), one can consider~\cite{Paper3,Bernstein2019feedback,dixit2019online}:
\begin{align}
\label{eq:timevarying_sequence}
\sigma_t := \|\bx_{t}^* - \bx_{t-1}^*\| , \hspace{.5cm} \Sigma_t := \sum_{\tau = 1}^{t} \sigma_\tau \, ,
\end{align}
with $\Sigma_t$  referred to as ``path length'' or ``cumulative drifting.'' The metric~\eqref{eq:timevarying_sequence} can be  utilized also when $f_t$ is a function of random parameters drawn from  a time-varying distribution; see, e.g.,~\cite{Cao2019}. When the cost is convex but not strongly convex,~\eqref{eq:timevarying_sequence} refers to a non-unique path and its respective length; however, an alternative measure that resolves this ambiguity involves a notion of worst-case path length~\cite{besbes2015non}:
\begin{align}
\label{eq:timevarying_sequence_worst_case} \Sigma_t^{\textrm{m}} := \max_{\{\bx_{\tau}^* \in \cX_\tau^*, \tau = 1, \ldots, t\}} \sum_{\tau = 1}^{t} \|\bx_{\tau}^* - \bx_{\tau-1}^*\| \, ,
\end{align}
where $\cX_\tau^*$ denotes the set of optimal solutions at time $\tau$. 

Additional metrics have been considered to capture the temporal variability of the underlying problem; for example, a variant of~\eqref{eq:timevarying_sequence} involving a dynamical model  is proposed in~\cite{hall2015online}. As another example, the metric $\Sigma_t^{\textrm{c}}(\bw) := \sum_{\tau = 1}^{t} \|\bw -
\bx_\tau\|^2$, where $\bw$ is a suitable ``comparator''; for example, $\bw$ could be the center of the solution trajectory $\bw = (1/\tau) \sum_{i = 1}^\tau \bx_i^*$ to capture the diameter of the minimizer sequence. For completeness, we also mention that, assuming that the constraint set $\cX_t$ is static (that is, $\cX_t = \cX$ for all $t$), then additional metrics are~\cite{besbes2015non,Jadbabaie2015}  
\begin{align}
\label{eq:timevarying_cost}
\Sigma_t^{\textrm{f}} := \sum_{\tau = 1}^{t} \sup_{\bx \in \cX} |f_\tau(\bx) -
f_{\tau-1}(\bx)|, \textrm{~and~} \,\, \Sigma_t^{\textrm{g}} := \sum_{\tau = 1}^{t} \sup_{\bx \in \cX} |\nabla h_\tau(\bx) -
\nabla h_{\tau-1}(\bx)|,
\end{align}
where, however, compactness of $\cX$ is needed to ensure finite value of $\Sigma_t^{\textrm{f}}$.

\subsection{Strongly convex time-varying functions: Linear convergence} 

We start by considering the case where the function $h_t$ in the cost of~\eqref{eq:time_varying_problem}  is $\mu_t$-strongly convex\footnote{If $g_t$ is $\mu_t$-strongly convex, we can assume $h_t$ is $\mu_t$-strongly convex without loss of generality since we can add $q(\bx_t):=\mu_t/2\|\bx_t\|^2$ to $h_t$ and subtract $q$ from $g_t$, and calculation of gradients and proximity operators of the new functions are given by standard formulas~\cite{CombettesBook2}.} and $L_t$-smooth, with $\mu_t \leq L_t$. Strong convexity implies that the sequence of optimal solutions $\{\bx_{t}^*, t \in \cT\}$ is unique; therefore, a pertinent performance assessment  involves the  analysis of the  \emph{tracking error} sequence $\{\|\bx_t - \bx_{t}^*\|_2, t \in \cT\}$~\cite{Bernstein2019feedback,Paper3}. To this end, one can
obtain
a slightly modified notion of $Q$-linear convergence of the form~\cite{Bernstein2019feedback,Paper3}:
\begin{align}
\label{eq:Qconvergence}
\|\bx_t - \bx_{t}^*\|_2 \leq q_t \|\bx_{t-1} - \bx_{t-1}^
*\|_2 + \cB_t(e_t,  \sigma_t)
\end{align}
\noindent for some $t \geq t_0$ (that is, without considering the transient behavior of the algorithm), where $q_t \in (0, 1)$ is the contraction coefficient, and $\cB_t: \mathbb{R}^2 \rightarrow \mathbb{R}$ is a function of the drifting $\sigma_t$ and the gradient error $e_t$. Assuming that there exists a scalar $\cB(e,  \sigma)$ that upper bounds the elements of the sequence $\{\cB_t(e_t,  \sigma_t)\}$, where  $\sigma < + \infty$ and $e < + \infty$ are bounds for  $\sigma_t$ and $e_t$ for all $t$, respectively,~\eqref{eq:Qconvergence} naturally leads to the following asymptotic result: 
\begin{align}
\label{eq:Qconvergence_asymptotic}
\limsup_{t \rightarrow + \infty} 
\|\bx_t - \bx_{t}^*\|_2
\leq \frac{\cB(e,  \sigma)}{1 - q}
\end{align}
where $q := \sup\{q_t\}$, which bounds the maximum tracking error. 

Concrete expressions for~\eqref{eq:Qconvergence_asymptotic} are provided in the following  two examples. 

$\bullet$ \emph{{\bf Online proximal method and projected gradient method}.} For the  example~\eqref{eq:forward-backward-v}, one has that  $q_t = \max\{|1 - \alpha_t \mu_t|,|1 - \alpha_t L_t|\}$~\cite{Be17}; therefore,~\eqref{eq:Qconvergence} is  a contraction (i.e., $q_t < 1$) when $\alpha_t < 2/L_t$. Taking a constant step-size $\alpha < 2/\sup\{L_t\}$, one has that  the bound~\eqref{eq:Qconvergence_asymptotic} is given by    

\vspace{-0.8cm}

\begin{align}
\label{eq:Qconvergence_asymptotic_fb}
\limsup_{t \rightarrow + \infty} \|\bx_t - \bx_{t}^*\|_2
\leq \frac{\alpha e + q \sigma}{1 - q} \, .
\end{align}
When $\sigma = 0$, one recovers results for the convergence of batch algorithms with errors in the gradient. Note also that this result allows for an approximate gradient ($e>0$) calculation.

We refer the reader to~\cite{necoara2019linear} for more results on linear convergence in the static setting when the cost function satisfies some relaxed strong convexity conditions.

$\bullet$ \emph{{\bf Online primal-dual method}.} We now turn the attention to primal-dual gradient methods based on the
Lagrangian function~\eqref{eq:reg-lagrangian}. The linear convergence results above are modified in this case to account for both primal and dual variables; that is, we consider the tracking error sequence $\{\|\bz_t - \bz_{t}^*\|_2, t \in \cT\}$, where $\bz_t := [\bx_t^\sfT, \blambda_t^\sfT]^\sfT$. In this case, the definition of the  drifting $\sigma_t$ is also modified accordingly as $\sigma_t^{\textrm{pd}} := \|\bz_{t}^* - \bz_{t-1}^*\|$. 

Still assuming that $f_t$ is strongly convex, the traditional (un-regularized) Lagrangian function is not strongly concave in $\blambda_t$ and thus no linear convergence results of the form~\eqref{eq:Qconvergence_asymptotic} may be possible~\cite{Koshal11,Du2019} (we notice also that in~\cite{Du2019} $f_t$ is not strongly convex, but a special structure of the regularized Lagrangian is assumed). When $r > 0$, then the regularized Lagrangian in~\eqref{eq:reg-lagrangian} is a strongly-convex strongly-concave function, and linear convergence results becomes available for both static~\cite{Koshal11,Du2019} and time-varying optimization problems~\cite{Bernstein2019feedback}; the price to pay, though, is tracking of the unique saddle-point of the \emph{regularized} Lagrangian function, which does not coincide in general with an optimal primal-dual pair of~\eqref{eq:time_varying_problem}; see, for example, the results in~\cite{Koshal11}. 
Sacrificing optimality for convergence 
is often appropriate in a time-varying setting, 
since if the regularization error
is small compared to the drift $\sigma_t$ and gradient error $e_t$, then a time-varying regularized algorithm would achieve very similar performance as a non-regularized one, with the added value of linear convergence. 

\newcommand{\mur}{\underbar{$\mu$}_t}
Letting $\bz_{t}^*$ be the saddle-point of the regularized Lagrangian when $r > 0$, and focusing first on the case where $g_t = 0$, the primal-dual  
operator 
$\mathcal{A}_t: (\bx,\blambda) \mapsto (\nabla h_t(\bx), -\bc_t(\bx) + r\blambda)$
is strongly monotone with strong-monotonicity constant $\mur = \min\{\mu_t, r\}$
(i.e., $\mathcal{A}_t-\mur I$, with $I$ the identity, is monotone), and Lipschitz continuous with a given constant $L_t^{\textrm{pd}}$  whenever $\bc_t$ is convex and  continuously  differentiable, and with a  Lipschitz  continuous  gradient (a precise derivation of $L_t^{\textrm{pd}}$ is available in e.g.,~\cite{Koshal11,Bernstein2019feedback}). Under these premises,  
$\limsup_{t \rightarrow + \infty} \|\bz_t - \bz_t^*\|_2$
can be bounded as in~\eqref{eq:Qconvergence_asymptotic_fb} with $\sigma_t$ replaced by $\sigma_t^{\textrm{pd}}$, $e$ an upper bound on the norm of the error in the computation of the gradient of both primal and dual steps~\cite{Bernstein2019feedback}, and $q_t = [1 - 2 \alpha \mur 
+ (L_t^{\textrm{pd}})^2]^\frac{1}{2} $.  Clearly, $q_t< 1$ if $\alpha$ is selected as $\alpha < 2 \mur 
/(L_t^{\textrm{pd}})^2$. These results hold also for the case with non-differentiable function $g_t$ and the proximal operator in the primal step, as shown in Eq.~\eqref{eq:primal-dual}.

\subsection{Dynamic regret} 

The dynamic regret is defined as~\cite{hall2015online,Jadbabaie2015,besbes2015non,chen2018bandit}: 
\begin{align}
\label{eq:dynamic_regret}
\textrm{Reg}_t := \sum_{\tau = 1}^t f_t(\bx_t) - \sum_{\tau = 1}^t  f_t^*
\end{align}
where we recall that $f_t^* = f_t(\bx_t^*)$ is in our case the optimal value function obtained by utilizing a batch algorithm. 
For strongly convex functions, boundedness of  $\|\bx_t - \bx_t^*\|$ implies boundedness of the instantaneous regret $f_t(\bx_t) -  f_t(\bx_t^*)$ when 
 $f_t$ is  Lipschitz continuous uniformly in $t$ (e.g., 
the (sub-)gradient of $f_t$ is bounded over the set $\cX_t$)~\cite{dixit2019online}; therefore, a recursive application of~\eqref{eq:Qconvergence} gives a bound on the dynamic regret. For the algorithm~\eqref{eq:forward-backward-v} and its projected gradient counterpart, it holds that $\textrm{Reg}_t = \mathcal{O}(1 + E_t + \Sigma_t)$~\cite{dixit2019online}. If the path length and the cumulative error are linear in  $t$,   no sublinear regret is possible as is confirmed 
by the lower bounds provided in, e.g.,~\cite{besbes2015non}. 

When the cost function is not strongly convex and the relaxed conditions for linear convergence explained in e.g, ~\cite{necoara2019linear} (for example, a quadratic functional
growth condition) are not satisfied, the dynamic regret then becomes a key performance metric. This is also the case for constrained problems, when the un-regularized Lagrangian is utilized to design the algorithm (that is, when i.e., $r = 0$ in~\eqref{eq:reg-lagrangian}), since even if the cost is strongly convex, the primal-dual operator is monotone (but not strongly monotone) if $r = 0$. A number of results in the literature are available for step-sizes $\alpha_t \propto 1/t$; here, on the other hand, we recall that we consider regret bounds for algorithms with a constant step-size
since they more naturally fit the time-varying setting.
 As an example, consider a smooth cost and set $\alpha = 1/L$, with $L_t \leq L$ for all $t$; then, for 
 an arbitrary ``comparator'' $\bw_t$ (i.e., reference for the performance analysis), 
 a bound for the  dynamic regret amounts to: 
\begin{align}
\label{eq:regret_bound}
    \textrm{Reg}_t 
    \leq \frac{L}{2}\left(\|\bx_0 - \bw_t\|^2+ \Sigma_t^{\textrm{c}}(\bw_t)\right)+\frac{1}{2L} E_t^\prime \, .
\end{align}
For example, if  $\bw_t = (1/\tau) \sum_{i = 1}^\tau \bx_i^*$, then $\Sigma_t^{\textrm{c}}(\bw_t)$ is the diameter of the minimizing sequence. Other bounds could be derived for approximate first-order information by extending the results of~\cite{besbes2015non,hall2015online}; they are close in spirit to~\eqref{eq:regret_bound}, and they imply  that no sublinear regret is achievable if the metric utilized to capture the time-variability of the problem grows  linearly. For completeness, we refer the reader to the lower bounds in, e.g.,~\cite{besbes2015non,li2018using}, and the bounds for primal-dual methods designed based on the standard Lagrangian function  in, e.g.,~\cite{chen2018bandit} and~\cite{Bernstein2019feedback}.

%%%%%%%%%%%%%%%%%%%%%%%%%%%%%%%%%%%
\section{Distributed Computation for information streams}
\label{sec:distributed}
%%%%%%%%%%%%%%%%%%%%%%%%%%%%%%%%%%%

Another key aspect of data streams is that they can be  distributed across different locations and sources.  With the increasing sheer amount of data, possibly coupled with privacy concerns, distributed computation plays a crucial role to ensure  that the data points are processed as close as possible to where they are generated.  We focus  on two key features and challenges in distributed time-varying optimization  with distributed information streams: \emph{(i)} step-size conditions; and,  \emph{(ii)} asynchronicity of the updates. Other aspects (e.g., communication vs. convergence, quantization, federated architectures) are also important in distributed optimization, in par with standard static processing; we will comment on these aspects in Section~\ref{sec:outlook} as part of the outlook.

Step-size conditions and synchronicity of the updates are two key differentiators between traditional static and online distributed optimization; if not handled properly, they may jeopardize performance and even convergence of standard distributed algorithms when applied to information streams. Take, for example,   decentralized (sub)-gradient descent (DGD): unless the step-size vanishes, convergence to the optimizer is not assured. On the other hand, as we discussed, if the step-size vanishes, then tracking of a time-varying optimizer becomes challenging. When considering cost functions that change over time, synchronicity becomes an even more important aspect; in principle, nodes at different locations are required to sample cost functions at the same instant, otherwise we would be solving  problems that pertain to different time instances at different nodes (jeopardizing performance at best, convergence at worst)\footnote{We focus on DGD, instead of distributed methods with possibly superior performance, for its simplicity, its chronological precedence, and because a number of current methods still employ it in some ways at their core. }.

\vspace{-.2cm}

\subsection{Example of time-varying consensus problem} 

To outline the ideas, consider a simplified version of~\eqref{eq:time_varying_problem} for a prototypical consensus problem: 

\vspace{-.8cm}

\begin{align}
\label{eq:time_varying_problem_distributed}
h_t^* := \min_{\bx_t \in \mathbb{R}^n} \sum_{i=1}^N h_{i,t}(\bx_t) \, .
\end{align}
Consider
$N$ spatially distributed nodes labeled as  $i = 1, \dots, N$, each one with a private cost $h_{i,t}(\cdot)$, which for simplicity of exposition we consider $L_t$-smooth and strongly convex. The nodes can communicate via a communication graph $\mathcal{G}$
and we will be looking at algorithms that would allow the nodes to agree on an optimizer $\bx^*_t$ of~\eqref{eq:time_varying_problem_distributed} at any time $t$, while communicating only with their neighboring nodes.  In the static setting, where the cost functions $h_{i,t}$ do not change over time, many algorithms have been developed to solve~\eqref{eq:time_varying_problem_distributed}~\cite{boyd2011distributed,dimakis2010gossip}, such as gradient tracking, exact first-order algorithm (EXTRA), dual averaging, dual decomposition, and ADMM (see, e.g.,~\cite{Shi2015} and references therein) in addition to DGD. We emphasize that the convergence claims of~\cite{boyd2011distributed,dimakis2010gossip,Shi2015} are for static optimization; a goal of this section is to highlight challenges in the design and analysis of distributed algorithms when moving from static optimization to time-varying optimization. Because of space limitations, we refer the reader to the work~\cite{Maros2019} for a comprehensive literature review on several aspects on time-varying distributed optimization, as well as~\cite{Hosseini2016,Akbari2017,Shahrampour2018} for examples of online and time-varying algorithms over networks.

DGD involves copies of the  variable $\bx$ to each node, denoted as $\by_i$, and generates the sequence
\vspace{-.3cm}
\begin{equation}\label{dgd}
    \by_{i, t} = \sum_{i \sim j} w_{ij} \by_{j, t-1} - \alpha_t \bv_{i,t}, \quad t = 1, 2, \dots
\end{equation}
where $i \sim j$ means that the sum is carried over all the neighbors of $i$ and $i$ itself, $w_{ij}$ are non-negative weights (which are often chosen as the relative degree between the nodes), $\alpha_t$ is the step-size, while $\bv_{i,t}$ is in this case the gradient of $h_{i,t}(\cdot)$ at $\by_{i, t-1}$, i.e., $\bv_{i,t} = \nabla_{\bx} h_{i,t}(\by_{i, t-1})$, or a proxy of the gradient as discussed in Section~\ref{sec:algorithms}. 
In the static setting, 
even in the strongly convex and $L_t$-smooth setting, the sequence $\{\by_{i, t}\}$ can be proven to converge, in the sense that $\lim_{t\to\infty} \|\by_{i, t} - \bx^*\|$, only when the step-size is vanishing: $\sum_{t} \alpha_t = \infty, \sum_{t} \alpha^2_t \leq \infty$ (under some extra but standard and mild assumptions on the communication graph, e.g., connectedness).
When the step-size is constant, then convergence is achieved only within a ball around of the optimizer. To understand this result, stack the variables $\by_{i, t}$ in a vector $\by_{t}$ and rewrite the recursion~\eqref{dgd} as

\begin{equation}\label{dgd-v}
    \by_{t} = \bW \by_{t-1} - \alpha_t \bv_{t}, \quad t = 1, 2, \dots
\end{equation}

\noindent where now $\bW$ is the matrix that contains the weights $w_{ij}$, while $\bv_{t}$ is the vector that contains all the local gradients. In particular, the matrix $\bI-\bW$ has maximum eigenvalue equal to $1$, with the corresponding eigenvector with all entries equal to $1$. Then, one can interpret~\eqref{dgd-v} as a standard gradient algorithm to solve the modified problem

\begin{align}
\label{eq:time_varying_problem_distributed-mod}
\min_{\by_t \in \mathbb{R}^n}  \sum_{i=1}^N h_{i,t}(\by_{i,t}) + \frac{1}{2 \alpha_t} \by_t^\sfT (\bI - \bW)  \by_t,
\end{align}
whose optimizer is different from~\eqref{eq:time_varying_problem_distributed} if $\alpha_t$ stays constant, showing that then exact convergence can never be achieved if  $\alpha_t$ is constant. 
It is also apparent that synchronicity must be enforced otherwise it is not clear what objective is being minimized; if the costs were sampled at different $t$'s at different nodes, the first term would read $\sum_{i=1}^N h_{i,t_i}(\by_{i,t})$ which is not the original objective.

\begin{figure}[h!]
    \centering
\includegraphics[width =1.0\textwidth]{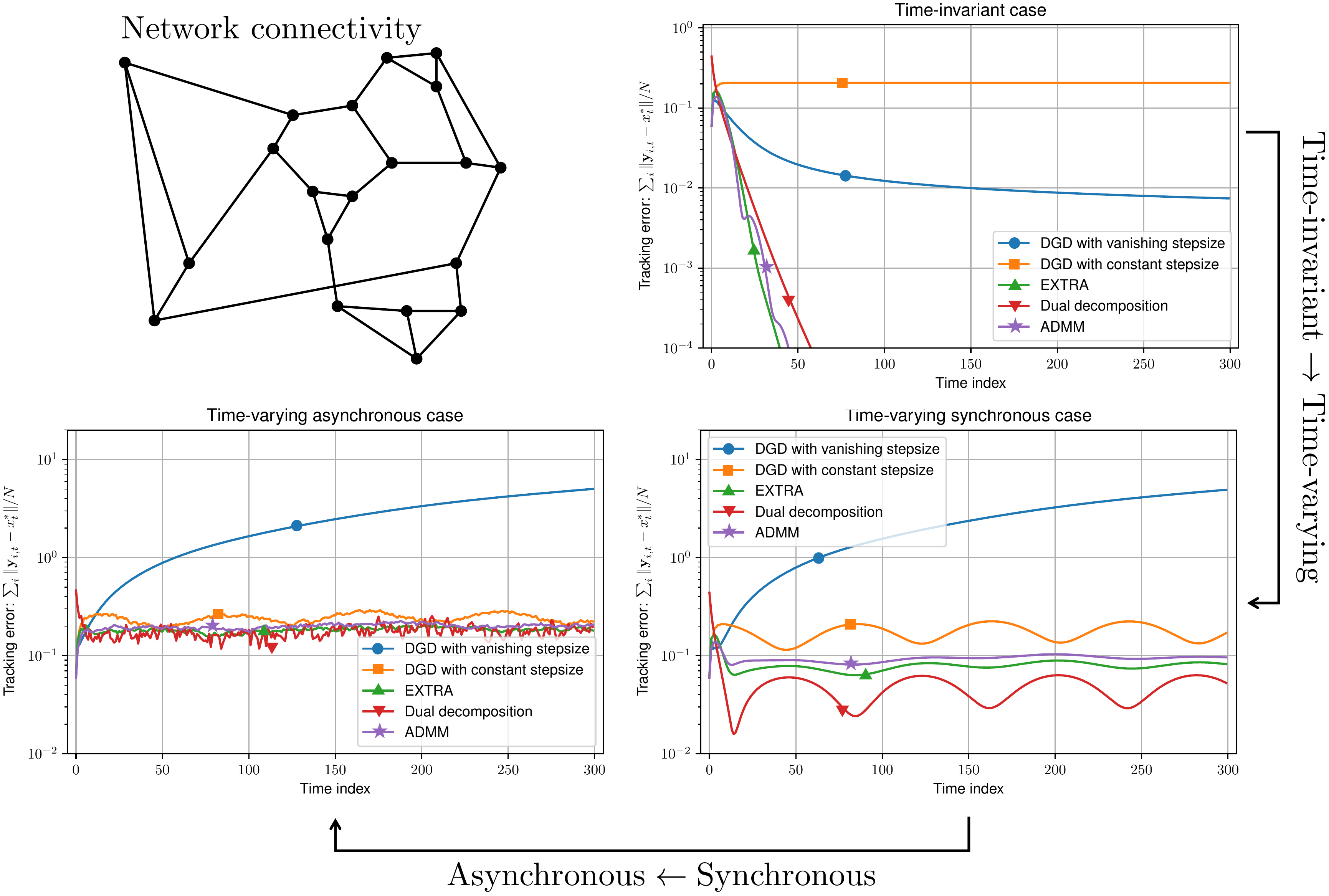}
    \caption{Numerical simulations for a time-varying optimization problem solved in a distributed way. \emph{Top left}: the communication graph consisting of 20 nodes; \emph{Top right}: tracking error in the time-invariant case;  \emph{Bottom right}: tracking error in the time-varying synchronous case; \emph{Bottom left}: tracking error in the time-varying asynchronous case.}
    \label{fig:distributed_sim_1}
\end{figure}

In Figure~\ref{fig:distributed_sim_1}, we illustrate the average tracking error $\sum_{i}\|\by_{i,t} - \bx_{t}^*\|/N$ for a time-varying problem defined over $N=20$ nodes (connectivity shown in upper left corner). 
The optimization problem is strongly convex and strongly smooth uniformly in time, and it has the form $\min_{x_t \in \mathbb{R}} \sum_{i=1}^N \frac{1}{2}\|x_t - A \cos(\omega t + \varphi_i - b t)\|^2 + \log(1 + \exp(x_t - a_i))$,
where $\{a_i\}$ and $\{\varphi_i\}$ are drawn 
i.i.d.~from uniform probability distribution of support $[-10,10]$ and $[0, 2\pi)$, respectively, while $A = 2.5,  \omega = \pi/80$. We study the performance of DGD with vanishing step-size ($\alpha = 1/t$), DGD with constant step-size, EXTRA, dual decomposition on the adjacency matrix of the graph, and distributed ADMM. Note that in this setting, the latter three methods have linear convergence in the static setting, and the latter two converge to an error bound in the time-varying setting~\cite{Ling14admm},~\cite{Jakubiec2013}. In the static setting (i.e., $\omega = b = 0$), EXTRA, dual decomposition, and ADMM maintain their theoretical promises and converge linearly. DGD with vanishing step-size converge slower, while DGD with constant step-size converge to an error bound. When we consider a time-varying setting, DGD with vanishing step-size diverges, while the other methods converge to an error bound as expected. Note that even if EXTRA has not been shown to converge in the time-varying setting, it is expected to do so, since it is linearly converging in the static setting. Having better performance in static setting 
does not clearly predict
better performance in the time-varying setting: for example, it seems that dual decomposition does much better in the time-varying scenario, while in the static setting it is worse than EXTRA and ADMM. 

Finally, 
the lower left corner of Fig.~\ref{fig:distributed_sim_1} illustrates
the case where we introduce asynchronicity in the sampling of the cost function (in this case nodes can sample functions asynchronously up to $10$ time instances in the past, meaning that each node has an $h_{i,t_i}$, with $t_i$ randomly generated and $0\leq t-t_i\leq 10$). The error is higher for all the algorithms, but it seems that DGD with constant step-size is the most robust. This is striking since DGD with constant step-size is the worst performing algorithm in the static setting, and shows once more that results in a static scenario cannot be easily translated into time-varying settings.

\section{Outlook}
\label{sec:outlook}

Streams of heterogeneous and spatially distributed data impose significant communication and computational  strains  on  existing  algorithmic  solutions. Deploying hardware with more powerful computational means is simply not a viable choice in many applications, and communication constraints still create severe bottlenecks in massively distributed settings. Time-varying optimization is rapidly emerging as an attractive solution. This article  emphasizes that we must revisit key design principles for batch optimization to enable an online processing of data without losing information or optimization capabilities.

\begin{figure}[h]
    \centering
    \includegraphics[width=\textwidth]{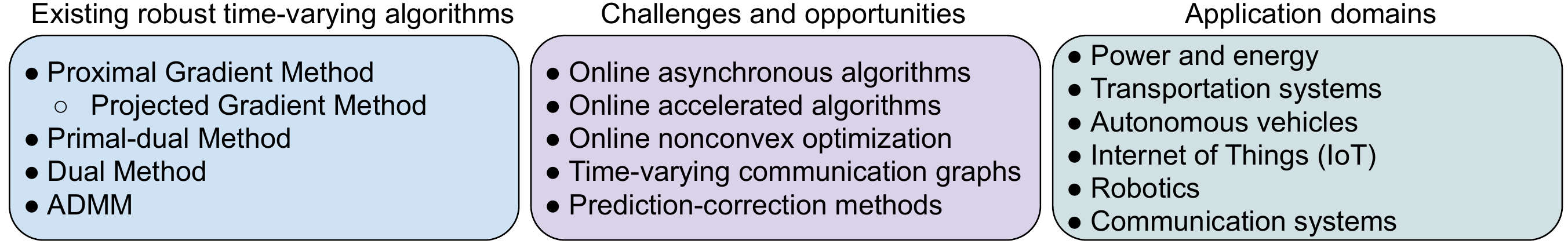}
    \caption{Outlook for analytical tools and application domains}
    \label{fig:outlook}
\end{figure}

Figure~\ref{fig:outlook} lists a number of key open questions in the time-varying optimization domain. For example, motivated by the representative results in Figure~\ref{fig:dynamic}, we  expect research efforts to explore the feasibility of accelerated methods in time-varying settings, along with efforts to characterize the performance of time-varying accelerated methods. Based on our discussion in Section~\ref{sec:distributed}, we also expect lines of research dealing with performance analysis of asynchronous time-varying algorithms, possibly operating with dynamic communications graphs. Lastly, while the present paper focuses on time-varying convex programs, we point out that a number of key research problems pertain to the development of  approximate online algorithms for time-varying nonconvex problems. This is driven by emerging applications such IoT and connected vehicles.  Figure~\ref{fig:outlook}  also lists a number of  application domains beyond the ML and SP areas.

%%%%%%%%%%%%%%%%%%%%%%%%%%%%%%%%%%%%%%%%%%%%%
\bibliographystyle{IEEEtran}
\bibliography{biblio.bib}

% Generated by IEEEtran.bst, version: 1.14 (2015/08/26)
\begin{thebibliography}{10}
\providecommand{\url}[1]{#1}
\csname url@samestyle\endcsname
\providecommand{\newblock}{\relax}
\providecommand{\bibinfo}[2]{#2}
\providecommand{\BIBentrySTDinterwordspacing}{\spaceskip=0pt\relax}
\providecommand{\BIBentryALTinterwordstretchfactor}{4}
\providecommand{\BIBentryALTinterwordspacing}{\spaceskip=\fontdimen2\font plus
\BIBentryALTinterwordstretchfactor\fontdimen3\font minus
  \fontdimen4\font\relax}
\providecommand{\BIBforeignlanguage}[2]{{%
\expandafter\ifx\csname l@#1\endcsname\relax
\typeout{** WARNING: IEEEtran.bst: No hyphenation pattern has been}%
\typeout{** loaded for the language `#1'. Using the pattern for}%
\typeout{** the default language instead.}%
\else
\language=\csname l@#1\endcsname
\fi
#2}}
\providecommand{\BIBdecl}{\relax}
\BIBdecl

\bibitem{friedman2001elements}
J.~Friedman, T.~Hastie, and R.~Tibshirani, \emph{The elements of statistical
  learning}.\hskip 1em plus 0.5em minus 0.4em\relax Springer, 2001.

\bibitem{cevher2014convex}
V.~Cevher, S.~Becker, and M.~Schmidt, ``Convex optimization for big data:
  Scalable, randomized, and parallel algorithms for big data analytics,''
  \emph{IEEE Signal Processing Magazine}, vol.~31, no.~5, pp. 32--43, 2014.

\bibitem{NocedalWright}
J.~Nocedal and S.~Wright, \emph{Numerical Optimization}, 2nd~ed.\hskip 1em plus
  0.5em minus 0.4em\relax Springer, 2006.

\bibitem{Be17}
A.~Beck, \emph{First-Order Methods in Optimization}.\hskip 1em plus 0.5em minus
  0.4em\relax {MOS}-{SIAM} Series on Optimization, 2017.

\bibitem{popkov2005gradient}
A.~Y. Popkov, ``Gradient methods for nonstationary unconstrained optimization
  problems,'' \emph{Automation and Remote Control}, vol.~66, no.~6, pp.
  883--891, 2005.

\bibitem{Zavala2010}
V.~M. Zavala and M.~Anitescu, ``{Real-Time Nonlinear Optimization as a
  Generalized Equation},'' \emph{SIAM Journal of Control and Optimization},
  vol.~48, no.~8, pp. 5444 -- 5467, 2010.

\bibitem{Bernstein2019feedback}
A.~{Bernstein}, E.~{Dall'Anese}, and A.~{Simonetto}, ``Online primal-dual
  methods with measurement feedback for time-varying convex optimization,''
  \emph{IEEE Trans. on Signal Processing}, vol.~67, no.~8, pp. 1978--1991,
  April 2019.

\bibitem{Paper3}
A.~Simonetto and E.~{Dall'Anese}, ``{Prediction-Correction Algorithms for
  Time-Varying Constrained Optimization},'' \emph{IEEE Transactions on Signal
  Processing}, vol.~65, no.~20, pp. 5481 -- 5494, 2017.

\bibitem{koller2018learning}
T.~Koller, F.~Berkenkamp, M.~Turchetta, and A.~Krause, ``Learning-based model
  predictive control for safe exploration,'' in \emph{2018 IEEE Conference on
  Decision and Control (CDC)}, 2018, pp. 6059--6066.

\bibitem{colombino2019online}
M.~Colombino, E.~Dall'Anese, and A.~Bernstein, ``Online optimization as a
  feedback controller: Stability and tracking,'' \emph{IEEE Transactions on
  Control of Network Systems}, 2019.

\bibitem{poveda2017framework}
J.~I. Poveda and A.~R. Teel, ``A framework for a class of hybrid extremum
  seeking controllers with dynamic inclusions,'' \emph{Automatica}, vol.~76,
  pp. 113--126, 2017.

\bibitem{dixit2019online}
R.~Dixit, A.~S. Bedi, R.~Tripathi, and K.~Rajawat, ``Online learning with
  inexact proximal online gradient descent algorithms,'' \emph{IEEE
  Transactions on Signal Processing}, vol.~67, no.~5, pp. 1338--1352, 2019.

\bibitem{elhamifar2013sparse}
E.~Elhamifar and R.~Vidal, ``Sparse subspace clustering: Algorithm, theory, and
  applications,'' \emph{IEEE Transactions on pattern analysis and machine
  intelligence}, vol.~35, no.~11, pp. 2765--2781, 2013.

\bibitem{Hajinezhad19}
D.~{Hajinezhad}, M.~{Hong}, and A.~{Garcia}, ``{ZONE}: Zeroth order nonconvex
  multi-agent optimization over networks,'' \emph{IEEE Transactions on
  Automatic Control}, 2019, early access.

\bibitem{chen2018bandit}
T.~Chen and G.~B. Giannakis, ``Bandit convex optimization for scalable and
  dynamic {IoT} management,'' \emph{IEEE Internet of Things Journal}, vol.~6,
  no.~1, pp. 1276--1286, 2018.

\bibitem{besbes2015non}
O.~Besbes, Y.~Gur, and A.~Zeevi, ``Non-stationary stochastic optimization,''
  \emph{Operations research}, vol.~63, no.~5, pp. 1227--1244, 2015.

\bibitem{onlineSaddle}
A.~Koppel, F.~Y. Jakubiec, and A.~Ribeiro, ``A saddle point algorithm for
  networked online convex optimization,'' \emph{IEEE Transactions on Signal
  Processing}, vol.~63, no.~19, pp. 5149--5164, Oct 2015.

\bibitem{Jadbabaie2015}
A.~Jadbabaie, A.~Rakhlin, S.~Shahrampour, and K.~Sridharan, ``{Online
  Optimization: Competing with Dynamic Comparators},'' in \emph{PMLR}, no.~38,
  2015, pp. 398 -- 406.

\bibitem{hall2015online}
E.~C. Hall and R.~M. Willett, ``Online convex optimization in dynamic
  environments,'' \emph{IEEE Journal of Selected Topics in Signal Processing},
  vol.~9, no.~4, pp. 647--662, 2015.

\bibitem{shalev2012online}
S.~Shalev-Shwartz \emph{et~al.}, ``Online learning and online convex
  optimization,'' \emph{Foundations and Trends{\textregistered} in Machine
  Learning}, vol.~4, no.~2, pp. 107--194, 2012.

\bibitem{hazan2016introduction}
E.~Hazan \emph{et~al.}, ``Introduction to online convex optimization,''
  \emph{Foundations and Trends{\textregistered} in Optimization}, vol.~2, no.
  3-4, pp. 157--325, 2016.

\bibitem{Dallanese2018feedback}
E.~Dall'Anese and A.~Simonetto, ``Optimal power flow pursuit,'' \emph{IEEE
  Transactions on Smart Grid}, vol.~9, no.~2, pp. 942--952, March 2018.

\bibitem{schmidt2011convergence}
M.~Schmidt, N.~L. Roux, and F.~R. Bach, ``Convergence rates of inexact
  proximal-gradient methods for convex optimization,'' in \emph{Advances in
  neural information processing systems}, 2011, pp. 1458--1466.

\bibitem{Koshal11}
J.~Koshal, A.~Nedi\'{c}, and U.~Y. Shanbhag, ``Multiuser optimization:
  Distributed algorithms and error analysis,'' \emph{{SIAM} J. on
  Optimization}, vol.~21, no.~3, pp. 1046--1081, 2011.

\bibitem{Ling14admm}
Q.~{Ling} and A.~{Ribeiro}, ``Decentralized dynamic optimization through the
  alternating direction method of multipliers,'' \emph{IEEE Transactions on
  Signal Processing}, vol.~62, no.~5, pp. 1185--1197, March 2014.

\bibitem{dimakis2010gossip}
A.~G. Dimakis, S.~Kar, J.~M. Moura, M.~G. Rabbat, and A.~Scaglione, ``Gossip
  algorithms for distributed signal processing,'' \emph{Proceedings of the
  IEEE}, vol.~98, no.~11, pp. 1847--1864, 2010.

\bibitem{boyd2011distributed}
S.~Boyd, N.~Parikh, E.~Chu, B.~Peleato, J.~Eckstein \emph{et~al.},
  ``Distributed optimization and statistical learning via the alternating
  direction method of multipliers,'' \emph{Foundations and
  Trends{\textregistered} in Machine learning}, 2011.

\bibitem{CombettesBook2}
H.~H. Bauschke and P.~L. Combettes, \emph{Convex Analysis and Monotone Operator
  Theory in Hilbert Spaces}, 2nd~ed.\hskip 1em plus 0.5em minus 0.4em\relax New
  York: Springer-Verlag, 2017.

\bibitem{halko2011finding}
N.~Halko, P.-G. Martinsson, and J.~A. Tropp, ``Finding structure with
  randomness: Probabilistic algorithms for constructing approximate matrix
  decompositions,'' \emph{SIAM Review}, vol.~53, no.~2, pp. 217--288, 2011.

\bibitem{trafficDBpaper2005}
A.~B. Chan and N.~Vasconcelos, ``Classification and retrieval of traffic video
  using auto-regressive stochastic processes,'' in \emph{IEEE Proceedings.
  Intelligent Vehicles Symposium, 2005.}\hskip 1em plus 0.5em minus 0.4em\relax
  IEEE, 2005, pp. 771--776.

\bibitem{Du2019}
S.~S. Du and W.~Hu, ``{Linear Convergence of the Primal-Dual Gradient Method
  for Convex-Concave Saddle Point Problems without Strong Convexity},'' in
  \emph{Proceedings of AISTATS}, 2019.

\bibitem{li2018using}
Y.~Li, G.~Qu, and N.~Li, ``Using predictions in online optimization with
  switching costs: A fast algorithm and a fundamental limit,'' in \emph{2018
  Annual American Control Conference (ACC)}, 2018, pp. 3008--3013.

\bibitem{Cao2019}
X.~{Cao}, J.~{Zhang}, and H.~V. {Poor}, ``On the time-varying distributions of
  online stochastic optimization,'' in \emph{2019 American Control Conference
  (ACC)}, July 2019, pp. 1494--1500.

\bibitem{necoara2019linear}
I.~Necoara, Y.~Nesterov, and F.~Glineur, ``Linear convergence of first order
  methods for non-strongly convex optimization,'' \emph{Mathematical
  Programming}, vol. 175, no. 1-2, pp. 69--107, 2019.

\bibitem{Shi2015}
W.~Shi, Q.~Ling, G.~Wu, and W.~Yin, ``{EXTRA: An Exact First-Order Algorithm
  for Decentralized COnsensus Optimization},'' \emph{SIAM Journal on
  Optimization}, vol.~25, no.~2, pp. 944 -- 966, 2015.

\bibitem{Maros2019}
M.~Maros, ``{Distributed Optimization in Time-Varying Environments},'' Ph.D.
  dissertation, KTH Stockholm, 2019.

\bibitem{Hosseini2016}
S.~{Hosseini}, A.~{Chapman}, and M.~{Mesbahi}, ``{Online Distributed Convex
  Optimization on Dynamic Networks},'' \emph{IEEE Transactions on Automatic
  Control}, vol.~61, no.~11, pp. 3545 -- 3550, 2016.

\bibitem{Akbari2017}
M.~{Akbari}, B.~{Gharesifard}, and T.~{Linder}, ``{Distributed Online Convex
  Optimization on Time-Varying Directed Graphs},'' \emph{IEEE Transactions on
  Control of Network Systems}, vol.~4, no.~3, pp. 417 -- 428, 2017.

\bibitem{Shahrampour2018}
S.~{Shahrampour} and A.~{Jadbabaie}, ``{Distributed Online Optimization in
  Dynamic Environments Using Mirror Descent},'' \emph{IEEE Transactions on
  Automatic Control}, vol.~63, no.~3, pp. 714 -- 725, 2018.

\bibitem{Jakubiec2013}
F.~Y. Jakubiec and A.~Ribeiro, ``{D-MAP: Distributed Maximum a Posteriori
  Probability Estimation of Dynamic Systems},'' \emph{IEEE Transactions on
  Signal Processing}, vol.~61, no.~2, pp. 450 -- 466, 2013.

\end{thebibliography}
%%%%%%%%%%%%%%%%%%%%%%%%%%%%%%%%%%%%%%%%%%%%%

%%%%%%%%%%%%%%%%%%%%%%%%%%%%%%%%%%%
\section{Authors}
\label{sec:biosketches}
%%%%%%%%%%%%%%%%%%%%%%%%%%%%%%%%%%%

{\small
\noindent \textbf{Emiliano Dall'Anese} (emiliano.dallanese@colorado.edu) is an assistant professor within the department of Electrical, Computer, and Energy Engineering at the University of Colorado Boulder. Previously, he was a postdoctoral associate at the University of Minnesota under the supervision of Prof. Georgios B. Giannakis, and a senior scientist at the National Renewable Energy Laboratory. He received the Ph.D. in Information Engineering from the University of Padova, Italy, in 2011. His research interests focus on optimization, decision systems, and signal processing, with applications to networked systems, hyperphysical systems, and power and energy systems.

\noindent \textbf{Andrea Simonetto} (andrea.simonetto@ibm.com) is a research staff member in the optimization and control group of IBM Research Ireland.
He received the Ph.D.\ degree in systems and control from Delft University of Technology, Delft, The Netherlands, in 2012. Previously, he was a postdoc at Delft University of Technology and at the Universit\`{e} catholique de Louvain, Louvain-la-Neuve, Belgium. He joined IBM Research in 2017. His interests span optimization, control, and signal processing, with applications in smart energy, smart transportation, personalized health, and quantum computing.

\noindent \textbf{Stephen Becker} (stephen.becker@colorado.edu) is an assistant professor in the Applied Mathematics department at the University of Colorado Boulder. Previously he was a Goldstine Postdoctoral fellow at IBM Research and a Fondation Scientifique et Mathematique de Paris postdoctoral fellow. He received his PhD in Applied Mathematics from Caltech under the supervision of Emmanuel Candes in 2011, and bachelors in math and physics from Wesleyan University.  His research interests are in optimization, machine learning, signal processing, imaging, inverse problems in quantum information, PDE-constrained optimization, tensor factorization and randomized numerical linear algebra.

\noindent \textbf{Liam Madden} (liam.madden@colorado.edu) is a Ph.D. student in the department of Applied Mathematics at the University of Colorado Boulder. He holds bachelor's degrees in Mechanical Engineering and Mathematics from California Polytechnic State University in San Luis Obispo, California.

}

\end{document}